\documentclass[eqsecnum,showpacs,showkeys,floats,aps,nofootinbib,preprint]{revtex4}
\usepackage{bm,amsfonts, mathtools}
\usepackage{graphicx,color, wasysym}
\usepackage{textcomp}
\usepackage{amsmath,amssymb,latexsym,epsfig}

\def\de{\mathrm{d}}

\def\nn{\nonumber}
\begin{document}

\makeatletter

\title{Symbolic calculus and integrals of Laguerre polynomials}

\author{D. Babusci}
\email{danilo.babusci@lnf.infn.it}
\affiliation{INFN - Laboratori Nazionali di Frascati, via E. Fermi, 40, IT 00044 Frascati (Roma), Italy}

\author{G. Dattoli}
\email{dattoli@frascati.enea.it}
\affiliation{ENEA - Centro Ricerche Frascati, via E. Fermi, 45, IT 00044 Frascati (Roma), Italy}

\author{K.~G\'{o}rska}
\email{kasia_gorska@o2.pl}
\affiliation{Instituto de F\'{\i}sica, Universidade de S\~{a}o Paulo, P.O.Box 66318, BR 05315-970 S\~{a}o Paulo, SP, Brasil}
\affiliation{H. Niewodnicza\'{n}ski Institute of Nuclear Physics, Polish Academy of Sciences, ul.Eljasza-Radzikowskiego 152, 
PL 31342 Krak\'{o}w, Poland}

\begin{abstract}
An umbral type formalism is used to derive integrals involving products of Laguerre polynomials and other special functions.
\end{abstract}

\maketitle


\section{Introduction}\label{s:intro}
It has been recently pointed out \cite{Anna} that, albeit integrals containing products of Gaussian functions and of an arbitrary 
number of Hermite polynomials are available since long time \cite{Prud}, the same is not true for Laguerre polynomials. A first step 
aimed at filling this gap has been put forward in \cite{Anna},  and in this paper we give our contribution using a recently 
developed \cite{DP} symbolic method, which has proved its usefulness in the study of integrals involving Bessel functions. 

Following this symbolic approach, the two-variable Laguerre polynomials 
\begin{equation}
L_n (x, y) = n!\,\sum_{k = 0}^n (- 1)^k\,\frac{x^{\,k}\,y^{\,n - k}}{(n - k)!\,(k!)^2}
\end{equation}
can be written as
\begin{equation}
\label{e:lague}
L_n (x, y) = (y - \hat{c}\,x)^n\,\varphi_0
\end{equation}
where the following umbral notation has been used \cite{DP,Datto}
\begin{equation}
\label{e:umbra}
\hat{c}^{\,\gamma}\,\varphi_0 = \varphi_\gamma = \frac1{\Gamma (\gamma + 1)}\,.
\end{equation}
$\varphi_0$ is referred to as the \textit{polynomial vacuum} (the umbral operator $\hat{c}$ acting on the ``state" $\varphi_0$ 
``creates" the Laguerre polynomials).

Before closing this short introduction, we report a result that will be useful in what follows. As it has been done in Ref. \cite{BDGermano}, 
by using the generating function of the two-variable Hermite polynomials
\begin{equation}
\label{e:Herm}
\sum_{n = 0}^\infty \frac{t^{\,n}}{n!}\,H_n (x,y) = e^{\,x\,t + y\,t^2} \qquad\qquad 
H_n (x,y) = n!\,\sum_{k = 0}^{[n/2]} \frac{x^{\,n - 2\,k}\,y^{\,k}}{(n - 2\,k)!\,k!}\,,
\end{equation}
it is possible to show that
\begin{equation}
\label{e:Inab}
I_n (a, b; \alpha, \beta) = \int_{- \infty}^\infty \de x \,(a\,x + b)^n\,e^{- \alpha\,x^2 + \beta\,x} = 
\sqrt{\frac{\pi}{\alpha}}\,\exp\left(\frac{\beta^2}{4\,\alpha}\right)\,H_n \left(b + \frac{a\,\beta}{2\,\alpha}, \frac{a^2}{4\,\alpha}\right)\,.
\end{equation}
The use of the symbolic method along with integrals of the type (\ref{e:Inab}) will be the key note for the future development.

\section{A survey of integrals of Laguerre polynomials}\label{s:lague}
By taking into account the expression \eqref{e:lague} for the Laguerre polynomials, one can write
\begin{equation}
{\cal I}_n (u; \alpha) = \int_{- \infty}^\infty \de x \,L_n (x, u)\,e^{- \alpha\,x^2} = I_n (- \hat{c}, u; \alpha, 0)\,\varphi_0\,,
\end{equation}
and, treating $\hat{c}$ as an ordinary constant, it's easy to show that
\begin{equation}
{\cal I}_n (u; \alpha) = \sqrt{\frac{\pi}{\alpha}}\,Q_n^{(0)} \left(u, \frac1{4\,\alpha}\right)
\end{equation}
where we have introduced the polynomials
\begin{align}
\label{e:Qpoly}
Q_n^{(\nu)} (x, y) &=  { \hat{c}^{\,\nu}\, H_n (x, y\,\hat{c}^{\,2})\,\varphi_0} \nn \\
&= n!\,\sum_{k = 0}^{[n/2]} \frac{x^{\,n - 2\,k}\,y^{\,k}}{(n - 2\,k)!\,k!\,\Gamma (2\,k + \nu + 1)}\,,
\end{align}
whose properties will be discussed later in the paper. 

The form \eqref{e:lague} generalizes to the case of associated two-variable Laguerre polynomials as follows
\begin{equation}
L_n^{(\nu)} (x, y) = \frac{\Gamma (n + \nu + 1)}{n!}\,\hat{c}^{\,\nu}\,(y - \hat{c}\,x)^n\,\varphi_0
\end{equation}
and, proceeding as above, we get
\begin{equation}
{\cal I}_n^{(\nu)} (u; \alpha) = \int_{- \infty}^\infty \de x \,L^{(\nu)}_n (x, u)\,e^{- \alpha\,x^2} = 
\sqrt{\frac{\pi}{\alpha}}\,\frac{\Gamma (n + \nu + 1)}{n!}\,Q_n^{(\nu)} \left(u, \frac1{4\,\alpha}\right)\,.
\end{equation}

These results and the identity\footnote{It is essentially a Taylor series expansion, and can be proved by noting that
$$
H_n (x + a, y) = e^{a\,\partial_x}\,H_n (x, y) =  \sum_{k = 0}^\infty \frac{(a\,\partial_x)^k}{k!}\,H_n (x, y)
$$
and 
$$
\partial_x^k\,H_n (x,y) = \frac{n!}{(n - k)!}\,H_{n - k} (x, y).
$$}
\begin{equation}
H_n (x + a, y) = \sum_{k = 0}^n \binom{n}{k}\,a^{\,k}\,H_{n - k} (x,y)
\end{equation} 
allow us to obtain the following formula
\begin{align}
_a{\cal I}_n^{(\nu)} (u; \alpha) & = \int_{- \infty}^\infty \de x \,L^{(\nu)}_n (x + a, u)\,e^{- \alpha\,x^2} \nn \\
& = \sqrt{\frac{\pi}{\alpha}}\,\frac{\Gamma (n + \nu + 1)}{n!}\,\sum_{k = 0}^n \binom{n}{k}\,(-a)^k\,Q_{n - k}^{(k + \nu)} 
\left(u, \frac1{4\,\alpha}\right)\,.
\end{align}

The multi-index extension of the Hermite polynomials can be exploited to compute integrals involving products of 
powers of binomials and gaussians in a compact form. In Ref. \cite{BD4} it has been shown that 
\begin{equation}
\int_{- \infty}^\infty \de x\,(a\,x + b)^m\,(f\,x + g)^n\,e^{\,- \alpha\,x^2} =\sqrt{\frac{\pi}{\alpha}}\,H_{m,n} 
\left(b, \frac{a^2}{4\,\alpha}; g, \frac{f^2}{4\,\alpha} \Big{|} \frac{a\,f}{2\,\alpha} \right)\,,
\end{equation}
where 
\begin{equation}
H_{m,n}  (x, y; w, z | \tau) = m!\,n!\,\sum_{k = 0}^{\mathrm{min} (m,n)} \frac{\tau^{\,k}}{(m - k)!\,(n - k)!\,k!}\,
H_{m - k} (x, y)\,H_{n - k} (w, z)
\end{equation}
are the two-index Hermite polynomials, whose properties are discussed in detail in Refs. \cite{Appe} and, more recently, 
in \cite{Datto}. The use of the same argument as before yields the following result
\begin{align}
{\cal I}_{m, n}^{(\mu, \nu)} (u, v; \alpha) & = \int_{- \infty}^\infty \de x \,L_m^{(\mu)} (x, u)\,L_n^{(\nu)} (x, v)\,e^{- \alpha\,x^2} \nn \\ 
 &= \sqrt{\frac{\pi}{\alpha}}\,\frac{\Gamma (m + \mu + 1)}{m!}\,\frac{\Gamma (n + \nu + 1)}{n!}\,
Q_{m, n}^{(\mu, \nu)} \left(u, \frac1{4\,\alpha}; v, \frac1{4\,\alpha} \Big{|} \frac1{2\,\alpha}\right)
\end{align}
where 
\begin{equation}
Q_{m, n}^{(\mu, \nu)} (x, y; w, z | \tau) = m!\,n!\,\sum_{k = 0}^{\mathrm{min} (m,n)} \frac{\tau^{\,k}}{(m - k)!\,(n - k)!\,k!}\,
Q_{m - k}^{(k + \mu)} (x, y)\,Q_{n - k}^{(k + \nu)} (w, z)\,.
\end{equation}
(The case of the product of two Laguerre polynomials is obtained putting $\mu = \nu = 0$). The result for an arbitrary product of 
associated Laguerre polynomials can be easily derived from the previous identities. It involves multi-index $Q$-polynomials and 
it is not reported for the sake of brevity.

The method is flexible enough to be easily extended to the calculation of family of integrals containing products of Laguerre and 
Hermite polynomials. Taking into account Eq. (9) of Ref. \cite{BD4}, one obtains 
\begin{equation}
\int_{- \infty}^\infty \de x\,(a\,x + b)^m\,H_n (f\,x +  g, y)\,e^{- \alpha\,x^2} = \sqrt{\frac{\pi}{\alpha}}\,
H_{m, n} \left(b, \frac{a^2}{4\,\alpha}; g, y + \frac{f^2}{4\,\alpha} \Big{|} \frac{a\,f}{2\,\alpha} \right)\,, 
\end{equation}
and, thus, by using expression \eqref{e:lague} for the associated Laguerre polynomials, it's easy to show that 
\begin{equation}
\int_{- \infty}^\infty \de x\,L_m^{(\nu)} (x, y)\,H_n (f\,x +  g, z)\,e^{- \alpha\,x^2} = 
\sqrt{\frac{\pi}{\alpha}}\,\frac{\Gamma (n + \nu + 1)}{n!}\,
T_{m, n}^{(\nu)} \left(y, \frac1{4\,\alpha}; g, z+ \frac{f^2}{4\,\alpha} \Big{|} - \frac{f}{2\,\alpha}\right)
\end{equation}
where 
\begin{equation}
T_{m, n}^{(\nu)} (x, y; w, z | \tau) = m!\,n!\,\sum_{k = 0}^{\mathrm{min} (m,n)} \frac{\tau^{\,k}}{(m - k)!\,(n - k)!\,k!}\,
Q_{m - k}^{(k + \nu)} (x, y)\,H_{n - k} (w, z)\,.
\end{equation}

As a final example, let us consider now the following integral
\begin{equation}
B_n (y, \alpha) = \int_{- \infty}^\infty \de x \,L_n (x, y)\,J_0 (2\,\sqrt{x})\,e^{- \alpha\,x^2} 
\end{equation}
where $J_0 (x)$ is the first-order cylindrical Bessel function, that, according to Ref. \cite{Datto}, can be cast in the form
\begin{equation}
J_0 (x) = \exp\left\{- \hat{d}\,\left(\frac{x}2\right)^2\right\}\,\lambda_0 
\end{equation}
with the umbral operator $\hat{d}$ acting on the polynomial vacuum $\lambda_0$ in the same way that $\hat{c}$ acts on 
$\varphi_0$ (see Eq. \eqref{e:umbra}). This integral looks quite complicated, but using the expression \eqref{e:lague} for 
the Laguerre polynomials, and taking into account Eq. \eqref{e:Inab}, we can write
\begin{equation}
B_n (y, \alpha) = \sqrt{\frac{\pi}{\alpha}}\,H_n\left(y + \frac{\hat{c}\,\hat{d}}{2\,\alpha}, \frac{\hat{c}^{\,2}}{4\,\alpha}\right)\,
\exp\left\{\frac{\hat{d}^{\,2}}{4\,\alpha}\right\}\,\varphi_0\,\lambda_0\,.
\end{equation}
Since
\begin{equation}
\label{e:Wright}
\hat{d}^{\,\nu}\,e^{\hat{d}^2\,x}\,\lambda_0 = \sum_{k = 0}^\infty \frac{x^{\,k}}{k!\,\Gamma (2\,k + \nu + 1)} = W_\nu (x | 2)\,,
\end{equation}
we get 
\begin{equation}
B_n (y, \alpha) = \sqrt{\frac{\pi}{\alpha}}\,\sum_{k = 0}^{[n/2]}  \frac{y^{\,n - 2 k}\,(4\,\alpha)^{- k}}{k!}\,\,A_{n - 2 k} (y, \alpha)
\end{equation}
where 
\begin{equation}
A_{n - 2k} (u, v) = \sum_{p = 0}^{n - 2 k} \binom{n}{2\,k + p}\,\frac{(2\,u\,v)^{- p}}{p!}\,W_p \left(\frac1{4\,v} \Big{|} 2\right)\,.
\end{equation}
The function $W_\alpha (x | 2)$ is the second-order Bessel-Wright function \cite{Andrews}, that intervenes in the expression 
of the generating function of the $Q$-polynomials defined in Eq. \eqref{e:Qpoly}, i.e. 
\begin{equation}
\sum_{n = 0}^\infty \frac{t^n}{n!}\,Q_n^{(\alpha)} (x, y) = \hat{c}^{\,\alpha}\,e^{\,x\,t + y\,t^2\,{\hat{c}^2}}\,\varphi_0 = 
e^{\,x\,t}\,W_\alpha (y\,t^2 | 2)\,,
\end{equation}
where the expression for the generating function of the Hermite polynomials (see Eq. \eqref{e:Herm}), and 
the identity \eqref{e:Wright} have been used. This means that $Q$-polynomials belong to the Sheffer family and 
can be defined by the operational identity
\begin{equation}
Q_n^{(\alpha)} (x, y) = W_\alpha (y\,\partial_x^2 | 2)\,x^n\,.
\end{equation}
They cannot be considered orthogonal polynomials in a strict sense of the word, but pseudo-orthogonal according to the 
notion introduced in \cite{Datto2}. This aspect of the problem will be addressed in a forthcoming publication.
%
%
%
%

\end{document}